# Validation of Flow Models as New Simulation Approach for Parcel Handling in Bulk Mode

*M. Sc. Domenik Prims* [a]
*M. Sc. Jennifer Kötz* [b]
*Prof. Dr. rer. nat. habil. Simone Göttlich* [b]
*Prof. Dr.-Ing. André Katterfeld* [c]

[a] Siemens Postal, Parcel & Airport Logistics GmbH,
Research & Development, Systems, Architecture
[b] University of Mannheim, Department of Mathematics
[c] Otto von Guericke University Magdeburg,
Institute of Logistics and Material Handling Systems

**F**low Models based on hyperbolic partial differential equations (conservation laws) are a well-established approach for the material flow simulation. Until now, they have been exclusively used for the simulation of cylindrical cargo. This paper investigates the application of flow models on cubical cargo as a new simulation approach for parcel logistics. Selected flow model parameters are adapted to cover this new situation. The simulation results of this macroscopic model are compared with the microscopic Discrete Element Method (DEM) where the cargo is approximated by superquadrics. An experimental setup especially designed for the validation of the considered flow model for cubical cargo bulk flow is presented. Results of the flow model are analyzed and validated against the results of experiment and DEM for the test setting.

*[Keywords: discrete element method, DEM, parcel handling in bulk mode, superquadrics, particle simulation, conservation laws, numerical studies]*

## 1  INTRODUCTION

The global courier, express and parcel (CEP) market is continuing to grow by 15 to 20 percent every year, according to the latest studies. Until 2022, an increase of 5.2 percent per year is expected in Germany. So the volume will rise up to more than 4.3 billion parcels per year [BIEK18]. The yearly European volume is assumed to reach nearly between 8 and 10 million swap body unloading processes.

Suppliers of parcel handling machinery like Siemens Postal, Parcel & Airport Logistics GmbH develop bulk handling technology like the Rubus (fully-automatic) to unload parcels in bulk mode. After the unloading process, the singulation of the parcels takes place as next process step. The parcels are conveyed on conveyor belts and often singulated by the support of different types of diverter equipment.

Regarding engineering issues, an usual approach would be to simulate the parcel bulk handling process based on the Discrete Element Method (DEM) to consider the real physical interaction of parcels before a new system is built up [PK17]. However, with the existing modelling approaches, the calculation effort increases with the number of parcels.

This paper covers the validation of a macroscopic model for material flow prediction in comparison to conventional parcel simulation approaches such as DEM. The calculation results of macroscopic models are independent from the number of items and are therefore suitable for the computationally efficient simulation of parcels in bulk mode, where the number of parcels is high.

## 2  EXPERIMENTAL SETUP

The experimental setup is needed for the validation of the microscopic DEM model (see Section 3) and the macroscopic flow model (see Section 4) for the bulk handling.

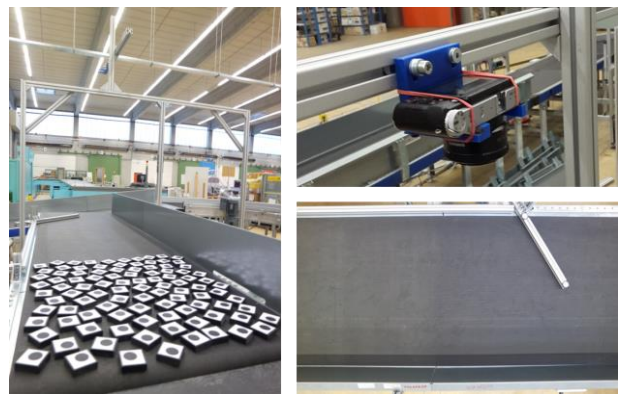

*Figure 1.   Siemens testing system for cargo bulk flow*

## 2.1 BULK FLOW TESTING SYSTEM

The experiments carried out by the R&D department of the Siemens Postal, Parcel & Airport Logistics GmbH cover the transport of cubical cargo on a conveyor belt (see Figure 1). Within the testing system, the cargo is redirected by a diverter. For the first analysis, two angles for the position of the diverter are evaluated: α = 45deg and α = 60deg (see Figure 2). The transport velocity of the conveyor belt is $v_T \approx 0.13$ m/s.

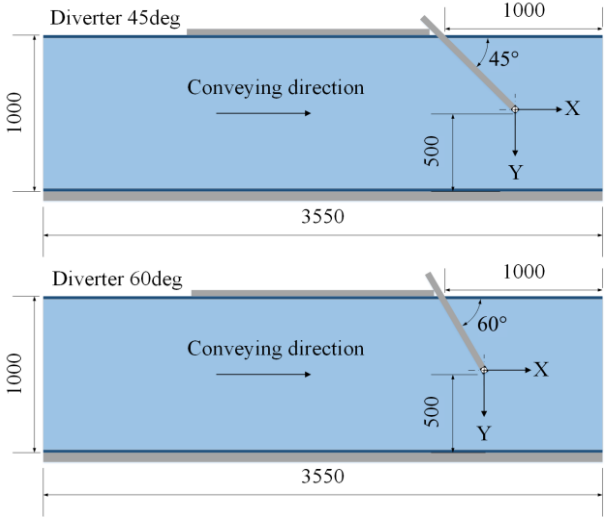

*Figure 2.   Layout for the experimental setup*

The end of the diverter marks a reference point which is used as starting point for the pattern recognition (see Figure 2) as well as threshold (x = 0) for the mass flow measurement of the microscopic and macroscopic simulation.

## 2.2 CARGO CHARACTERIZATION

For the first evaluation steps, the experiments are carried out with ideal rigid bodies because the real parcel rigidity is hard to describe. All the items consist of spruce wood (see Figure 3). Each item has l = 0.07m length, w = 0.07m width, h = 0.02m height and a weight of m = 0.0491kg.

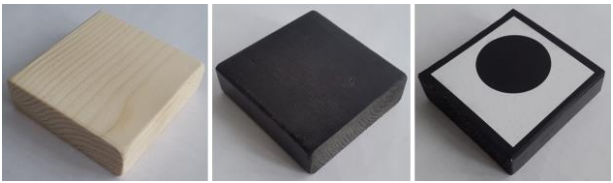

*Figure 3.   Cargo for experiments with pattern recognition*

The coefficient of static friction is determined for all contact surfaces of the experimental setup by a simple inclined plane test (see Table 1).

*Table 1.   Coefficients of static friction*

| Surface | Angle of static friction [°] | Coefficient of static friction [-] |
|---|---|---|
| Cargo | 30.6 | 0.59 |
| Diverter | 23.5 | 0.43 |
| Skirt boards | 30.4 | 0.59 |
| Conv. belt | 34.6 | 0.69 |

## 2.3 CARGO RECOGNITION

A special pattern was designed to recognize the cargo position and rotation during the transport (see Figure 4).

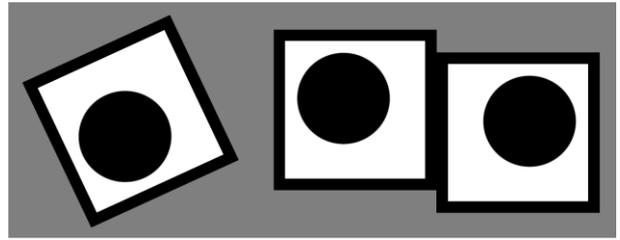

*Figure 4.   Cargo labelling for pattern recognition*

The black frame helps to distinguish between two adjacent items, while the black circle is used for tracking of the item rotation. Pattern recognition is also used for the determination of the conveyor speed in each experiment.

## 2.4 EXPERIMENTAL PROCEDURE

Each experiment involves N = 100 quadratic items. Figure 5 shows the initial positioning of items exemplary for the setup with a 45deg diverter.

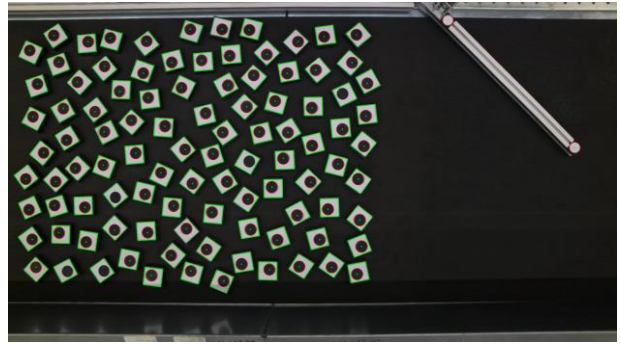

*Figure 5.   Initial positioning of items with item recognition frame (green color) and circle (red color)*

The experiments are repeated five times for each diverter angle (see Figure 6) and the initial cargo position and rotation is random. One representative experiment per diverter angle is selected for the comparison of microscopic and macroscopic results with experimental data:

- α = 45deg → experiment 0764,
- α = 60deg → experiment 0772.

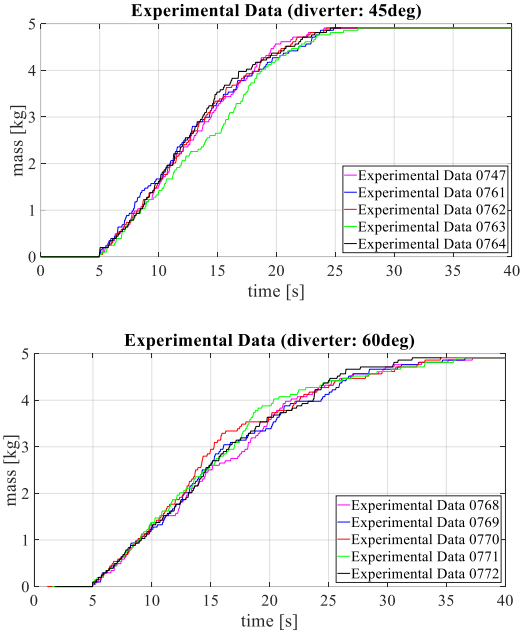

*Figure 6. Data of five experiments per diverter angle*

## 3 MATERIAL FLOW SIMULATION WITH DISCRETE ELEMENT METHOD

All microscopic simulations are carried out with LIGGGHTS® 3.8.0, an open source DEM tool.

### 3.1 SIMULATION MODEL

As the computation time of DEM software is increasing with the number of particles, it is even more critical to use the well-established multi-sphere approach because a high number of particles is needed to approximate the real parcel surface [SPPKP18]. In contrast, the superquadrics approach uses one particle per item which is computationally more efficient. As the superquadrics approach is available within LIGGGHTS® since 2017 [PPK17], it is applied for the cargo approximation within the DEM model.

The geometry of the experimental setup is reduced to the surfaces being in contact with the cargo during the transportation process (see Figure 7).

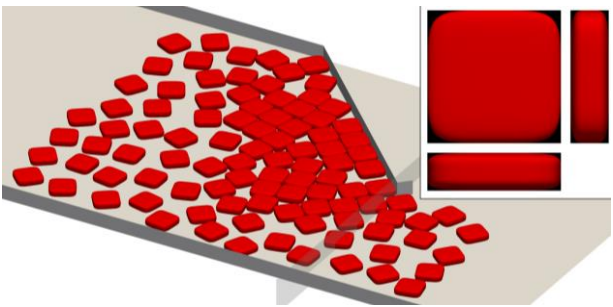

*Figure 7. DEM simulation of bulk flow with superquadrics, blockiness value $n_1 = n_2 = 5$ (see top right)*

The simulation parameters of the DEM contact model are shown in Table 2. Common guidelines say that a timestep below 15% of the Rayleigh time provides a sufficient accuracy. Therefore, the timestep of $\Delta t$ = 5e-5s has been chosen as it equals ~12% of Rayleigh time. The coefficient of restitution cor, the Poisson's ratio $\varepsilon$ and the Young's modulus E are based on simulation experience and do not significantly influence the DEM results. The calibration of the coefficients of sliding friction $\mu_s$ shows that the $\mu_s$ value needs to be set 50% lower than the coefficients of static friction (see Table 1) to match the real mass flow (see Figure 8).

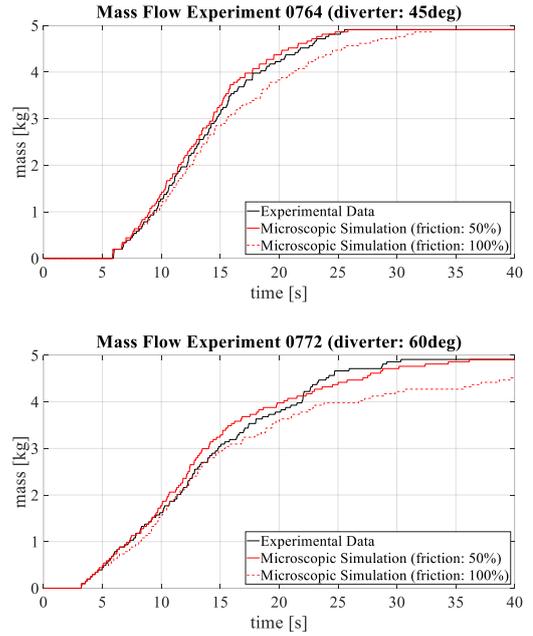

*Figure 8. DEM simulation vs. experiment 0764 ($\alpha$ = 45deg) and experiment 0772 ($\alpha$ = 60deg)*

However, the sliding friction value between the cargo and the conveyor belt is not reduced because of the influence of the ground curvature of the superquadrics on the cargo rotation during the transportation process.

*Table 2. Parameters of the DEM contact model*

| Parameter | Symbol | Unit | Value |
|---|---|---|---|
| Timestep | $\Delta t$ | s | 5e-5 |
| Coefficient of restitution | cor | - | 0.2 |
| Poisson's ratio | $\varepsilon$ | - | 0.3 |
| Young's modulus | E | N/m² | 1e+7 |
| Sliding friction | $\mu_s$ | - | see Table 1 (50%), 100% for conv. belt |
| Rolling friction | $\mu_r$ | - | 0.5 |
| Blockiness | $n_1$, $n_2$ | - | 5 |

There is no record that the rolling friction value $\mu_r$ has an influence on the bulk flow behavior. In general, computing time is increasing with the blockiness value. Here, the blockiness parameters $n_1$ and $n_2$ are low enough to keep computing time low and high enough to display the accurate cargo shape.

### 3.2 SIMULATION RESULTS

The mass flow is used as quantitative criteria to validate the DEM for bulk flow on conveyors with diverter equipment. The flow recording stops when the last one of the 100 items passes through the reference line (x = 0). In this moment, the total mass is m = 4.91kg. The simulation results match the experimental data for both diverter angles (see Figure 8) when the friction values are reduced by 50% – with one exception shown in Table 2.

During the data post processing, an uniform visualization style has been introduced to compare the experimental cargo behavior with microscopic data. The cargo position is represented by a rectangle and the orientation is shown by the circle inside the rectangle. Exemplary the comparison is shown for diverter angle α = 45deg (see Figure 12).

## 4 MACROSCOPIC ANALYSIS OF EXPERIMENTAL DATA

Macroscopic models based on conservation laws describe the material flow and treat items as a continuum. In this application, we use a two-dimensional hyperbolic partial differential equation to approximate the flow of cargo on the conveyor belt. The flow is characterized by macroscopic quantities such as the density and the velocity in space and time. In contrast to microscopic models like the DEM used before, the computing time of macroscopic models is independent of the number of items and therefore well suited for the experimental setups involving a huge number of items, i.e. in the magnitude of thousands.

### 4.1 ADAPTION OF THE MACROSCOPIC MODEL TO THE EXPERIMENTAL SETUP

We want to make use of the hyperbolic model for material flow [GHSV14] based on conservation laws to simulate the material flow of the experimental setup. The item density is defined as a function depending on space and time $\rho: \Omega \times \mathbb{R}^+ \to \mathbb{R}^+$, where $\Omega \subset \mathbb{R}^2$ is the two-dimensional state space with boundary $\delta\Omega$. We denote $\rho = \rho(x,t)$.

$$\partial_t \rho + \nabla \cdot \left(\rho \left(v^{dyn}(\rho) + v^{stat}(x)\right)\right) = 0,$$

$$v^{dyn}(\rho) = H(\rho - \rho_{max}) \cdot I(\rho),$$

$$I(\rho) = -\epsilon \frac{\nabla(\eta * \rho)}{\sqrt{1+||\nabla(\eta * \rho)||}},$$

$$\rho(x,0) = \rho_0(x).$$

Material flow is determined by the two components $v^{stat}$ and $v^{dyn}$. The static velocity field $v^{stat}(x)$ is induced by the geometry of the belt and its influence is determined according to the position x of the cargo. The dynamic vector field $v^{dyn}(\rho)$ includes the interaction between items weighted by the factor $\epsilon > 0$ and models the formation of congestion in front of an obstacle using the collision operator $I(\rho)$. The mollifier η which is used in the operator $I(\rho)$ is set as follows

$$\eta(x) = \frac{\sigma}{2\pi} e^{(-\frac{1}{2}\sigma\,||x||_2^2)}.$$

The Heaviside function H which activates the collision operator is approximated by

$$H(u) = \frac{1}{\pi}\arctan\bigl(h(u-1)\bigr) + \frac{1}{2}$$

where high values of h correspond to a better approximation of the Heaviside function.

### 4.2 TRANSFER TO THE EXPERIMENTAL SETUP

To numerically solve the conservation law, a grid in time and space with step sizes $\Delta_x$ and $\Delta_y$ for the state space and $\Delta_t$ for the time space is introduced. The solution of the macroscopic model is computed by the finite volume scheme with dimensional splitting described in [GHSV14]. A detailed overview of the parameters used is given in Table 3. The conveyor belt velocity for each setup is set to the experimentally measured velocity $v_T$.

*Table 3.   Parameters of the macroscopic model*

| Parameter | Symbol | Unit | Value |
|---:|:---:|:---:|:---:|
| Step size | $\Delta_x, \Delta_y$ | m | 0.01 |
| Timestep | $\Delta_t$ | s | 0.002 |
| Scaling factor | σ | - | 10,000 |
| Scaling factor | h | - | 50 |

The maximum density in the macroscopic model is

$$\rho_{max} = \frac{w}{l^2}$$

and is normalized to one for numerical computations. In the macroscopic simulation, each item is modeled as a square in the (x,y)-plane.

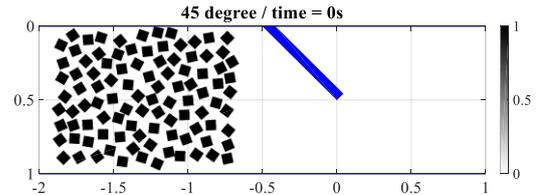

*Figure 9.   Initial density – macroscopic simulation*

The initial density corresponding to Figure 5 is depicted in Figure 9. The density takes its maximum in the inner surface of each item and is zero in between the items.

### 4.3 NUMERICAL RESULTS

To simulate the considered test setting, the space and time step sizes as well as the mollifier and the approximation of the Heaviside function are kept fixed for both diverter positions. Different parameters $\epsilon$ are applied for each scenario of the diverter and their results are contrasted below. The computational cost for the macroscopic model are constant for each of the diverter angles and depend only on the chosen discretization.

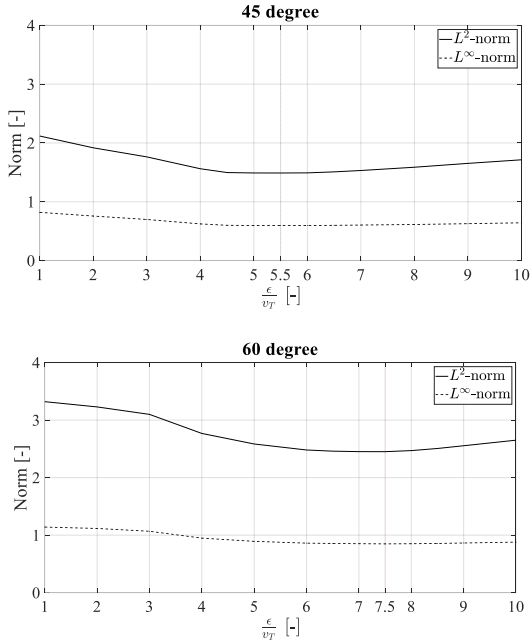

*Figure 10. $L^2$- and $L^\infty$-norm of macroscopic simulation data*

Figure 10 shows the error between the outflow computed with the simulation for different parameters $\epsilon$ and the outflow given by the experimental data (at x = 0). The error is measured in the discrete $L^2$- and $L^\infty$-norm. For both positions of the diverter, the error term decreases with increasing values for $\epsilon$ before reaching a minimum and increase afterwards. The minimum is reached for $\epsilon = 5.5 \cdot v_T$ and $\epsilon = 7.5 \cdot v_T$ for $\alpha = 45, 60$deg respectively. Too low or too high values of $\epsilon$ lead to a reduction in the accuracy of the macroscopic model. The simulation results show that the parameter $\epsilon$ should increase with increasing angles of the diverter to allow for a more precise description of the respective material flow.

Figure 13 (right column) shows the item density at different points in time during the simulation for $\alpha = 45$deg and $\epsilon = 5.5 \cdot v_T$. The figures show the distribution of the item density on the conveyor belt. Each density value is represented by one color. The black areas represent regions with maximum cargo density and the white areas represent regions without any items. At the beginning, the mass is transported on the conveyor belt with the belt velocity $v_T$. Afterwards, the mass starts to accumulate in front of the diverter (t = 5s). The region in front of the diverter gets more and more crowded and a region of maximal density evolves (t = 10s). Items are then redirected in the direction of the diverter and are able to leave the crowded region. At t = 25s most of the mass is able to leave the region in front of the diverter and the congestion in front of the obstacle is almost dissolved.

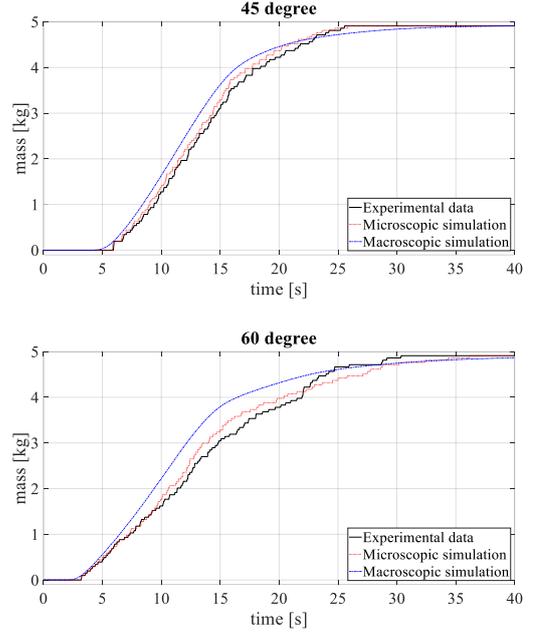

*Figure 11. Mass flow curve of the experiment, microscopic and macroscopic approach*

Figure 11 illustrates the outflow for $\alpha = 45$deg, $60$deg of the macroscopic simulation and the outflow given by the respective experiment. The parameter $\epsilon$ is set to $5.5 \cdot v_T$ and $7.5 \cdot v_T$ for $\alpha = 45, 60$deg. The outflow at time t is measured as the accumulated mass that has already passed the diverter at time t. Note that the mass flow of the simulation is a smooth curve in contrast to the experimental data, where mass flow is measured as the sum of the mass of all items whose center has exceeded the reference point. It is apparent from the figure that the mass flow curve fits the experimental data for $\alpha = 45$deg remarkably well. The congestion at the diverter between t = 5s and t = 20s dissolves faster in the macroscopic model than in the experiment. However, we observe that the mass flow curve of the macroscopic model is below the experimental curve at the end of the time frame.

The simulation results for the diverter positioned at $\alpha = 60$deg also fit the experimental data quite well. In contrast, due to the steepened angle of the diverter, we observe additional congestions and more crowded regions in front of the diverter. The simulation result is similar to the setup with the diverter positioned at 45deg. First, mass leaves the region faster in the macroscopic simulation than in the experiment, but after t = 25s mass leaves the region faster in the experiment. Due to the steeper angle, simulation results for the mass flow are not as close to the experimental data anymore, because the cubical items get stuck in front of the diverter, which is not included to the same extent in the macroscopic model.

The numerical results provide evidence that the macroscopic model can describe the experimental setup quite well. We achieve the most accurate results for the diverter positioned at 45deg. The error term in all considered norms is the lowest for this position of the diverter. The 60deg diverter still can be modeled well, while there are some losses in the accuracy of the description of the experiment. A possible explanation is that the macroscopic model was originally designed for cylindrical instead of cubical items. Therefore, the model cannot exactly portray the situation where items get stuck in front of the obstacle.

Note that we applied the numerical scheme described in [GHS14] to construct the solution to the conservation law. The numerical scheme uses a modified Roe flux, a mollifier and an approximation to the Heaviside function. Changes in the shape of the mollifier and the approximation of the Heaviside function do not significantly shift the mass flow curve closer to the experimental data.

## 5 COMPARISON BETWEEN MICROSCOPIC AND MACROSCOPIC RESULTS

Figure 12 and Figure 13 present the comparison of experimental data with the results of the microscopic DEM simulation and of the macroscopic simulation for the diverter angle α = 45deg. In contrast to visualization per item of the DEM simulation and the experimental data, the flow model shows a density distribution. Congestions in front of the diverter can already be observed at t = 5s. The black area, which marks the region of maximum density in the macroscopic model, mostly coincides with the region in the experiment and in the DEM simulation where the cargo is side by side. At t = 10s all of the images show additional congestions of items/ mass in front of the obstacle.

The mass flow of the microscopic and the macroscopic model matches shortly after t = 20s (see also Figure 11). At this point in time, the same cargo mass is left in front of the reference point (x = 0) in the DEM simulation as well as in the macroscopic simulation. Figure 14 shows that the cargo distribution in front of the diverter differs from the microscopic to the macroscopic simulation. One can see that the cargo positioning calculated with the DEM is closer to the experiment than the macroscopic simulation results. However, the two simulations predict the same mass flow at this time.

The last items pass the reference point at the end of the diverter at t = 25s in the experiment and in the DEM simulation. Compared to the same point in time in the macroscopic simulation, there is only a low density left in front of the diverter.

The runtimes of the macroscopic and microscopic simulation for the time horizon of 40 seconds and only 100 parts are compared on an Intel® Xeon® E5-2687W v3 CPU with 3.10 GHz. The computing time for a single simulation run is about 2 minutes in the microscopic model and about 37 minutes in the macroscopic model for the respective discretization given in Table 2 and Table 3. Our findings show that the flow model with the fine resolution given in Table 3 is not a suitable alternative to improve the computational efficiency of parcel bulk simulations when the item number is too low. For the presented test scenario, the computation time of the DEM model is only ~5% of the computation time of the flow model. But the flow model is able to portray the flow behavior of cubical cargo at material handling equipment and to estimate its mass flow even after redirection. First projections show that the flow model gets computationally more efficient when the item number exceeds 500.

Note that the resolution for the macroscopic model given in Table 3 is a very fine resolution. We already achieve convenient results for the mass flow for a coarser resolution ($\Delta_x = \Delta_y = 0.02m$, $\Delta_t = 0.004s$), where the macroscopic simulation time is about 4 minutes. Furthermore, the flow model becomes more competitive when smaller time step sizes in the microscopic simulation have to be considered. For a time step size of $\Delta t$ = 1e-5s, the computation time of the DEM amounts to ~8 minutes and therefore already exceeds the macroscopic computation time.

The selection of appropriate step sizes in each of the two models clearly depends on the level of accuracy and the desired runtime. We conclude that for the considered test setting, the macroscopic model is not a competitive alternative regarding the runtime because the item number is too low. However, the flow model is able to qualitatively represent the mass flow for the test data setting. The macroscopic model is a suitable alternative in test scenarios where qualitative statements on the mass flow of a high item number are of interest, because the runtime of the macroscopic model is independent of the number of items. Especially in the simulation of parcels in bulk mode, where the number of items is extremely high, the macroscopic model can serve as a substitute model to simulate the flow of parcels on a conveyor belt.

## 6 CONCLUSION

The microscopic Discrete Element Method and the macroscopic flow model were validated based on a real data test setting. A parameter configuration for the contact model of the DEM simulation was presented which leads to a very good correlation between simulation results and experimental data for different diverter angles. A macroscopic model which makes use of a partial differential equation to simulate the flow of items on a conveyor belt was adapted to the test setting and was validated against the experimental data and the results of the DEM. With the results of this study, we show that the flow model is suitable for the simulation of cubical cargo and especially for the simulation of a high item number in bulk mode. We will do further research on the macroscopic model to determine the most efficient combination of space and time step sizes. It is also proposed that further research is undertaken on the validation of flow models for the bulk flow of real parcels.

## FUNDING NOTE

The presented results in this paper are covered by the research projects "Simulation des Bewegungsverhaltens gefüllter Pakete und Ladungsträger im Pulk (SIMPPL)" and "Kombinierte Optimierung und Virtuelle Inbetriebnahme von materialflussintensiven Produktionssystemen mit Multiskalen-Netzwerk-Modellen (OptiPlant)". The projects No. KA1802/2-1 (SIMPPL) and G01920/7-1 (OptiPlant) are founded by the German Research Foundation (DFG).

───────────────────────────

**Domenik Prims, M. Sc.,** System Architect at the central Research & Development Department of the Siemens Postal, Parcel & Airport Logistics GmbH.

Domenik Prims was born in 1989 in Naumburg/S., Germany. Between 2008 and 2013 he studied Mechanical Engineering at the Otto von Guericke University Magdeburg. Between 2013 and 2015, he worked as head of department for material flow simulation at FAM Magdeburger Förderanlagen und Baumaschinen GmbH. In 2016, he worked as research associate at the Otto von Guericke University Magdeburg. Since 2016, he works for Siemens.

Address: Siemens Postal, Parcel & Airport Logistics GmbH, Lilienthalstraße 16/18, 78467 Konstanz, Germany

**Jennifer Kötz, M. Sc.,** PhD student at the University of Mannheim

Jennifer Kötz was born in 1993 in Cologne, Germany. She studied Business Mathematics at the University of Mannheim from 2012 to 2018. Since April 2018, she works as research assistant at the Scientific Computing Research Group.

**Simone Göttlich, Prof. Dr. rer. nat. habil.,** Professor at the School of Business Informatics and Mathematics at the University of Mannheim

Since 2011, Simone Göttlich directs the Scientific Computing Research Group at the University of Mannheim. Her research focuses on mathematical modeling, numerical simulation and optimization of dynamic processes. It includes various applications in the field of manufacturing systems, traffic flow, pedestrian dynamics and power grids.

Address: University of Mannheim, Department of Mathematics, 68131 Mannheim, Germany

**André Katterfeld, Prof. Dr.-Ing.,** Professor at the Otto von Guericke University Magdeburg

Since 2009, André Katterfeld is Chair of Material Handling Systems at the Institute of Logistics and Material Handling Systems of the Otto von Guericke University. One of his research topics is the calibration and application of DEM simulation in the field of material handling.

Address: Otto von Guericke University Magdeburg, ILM, Universitätsplatz 2, 39106 Magdeburg, Germany

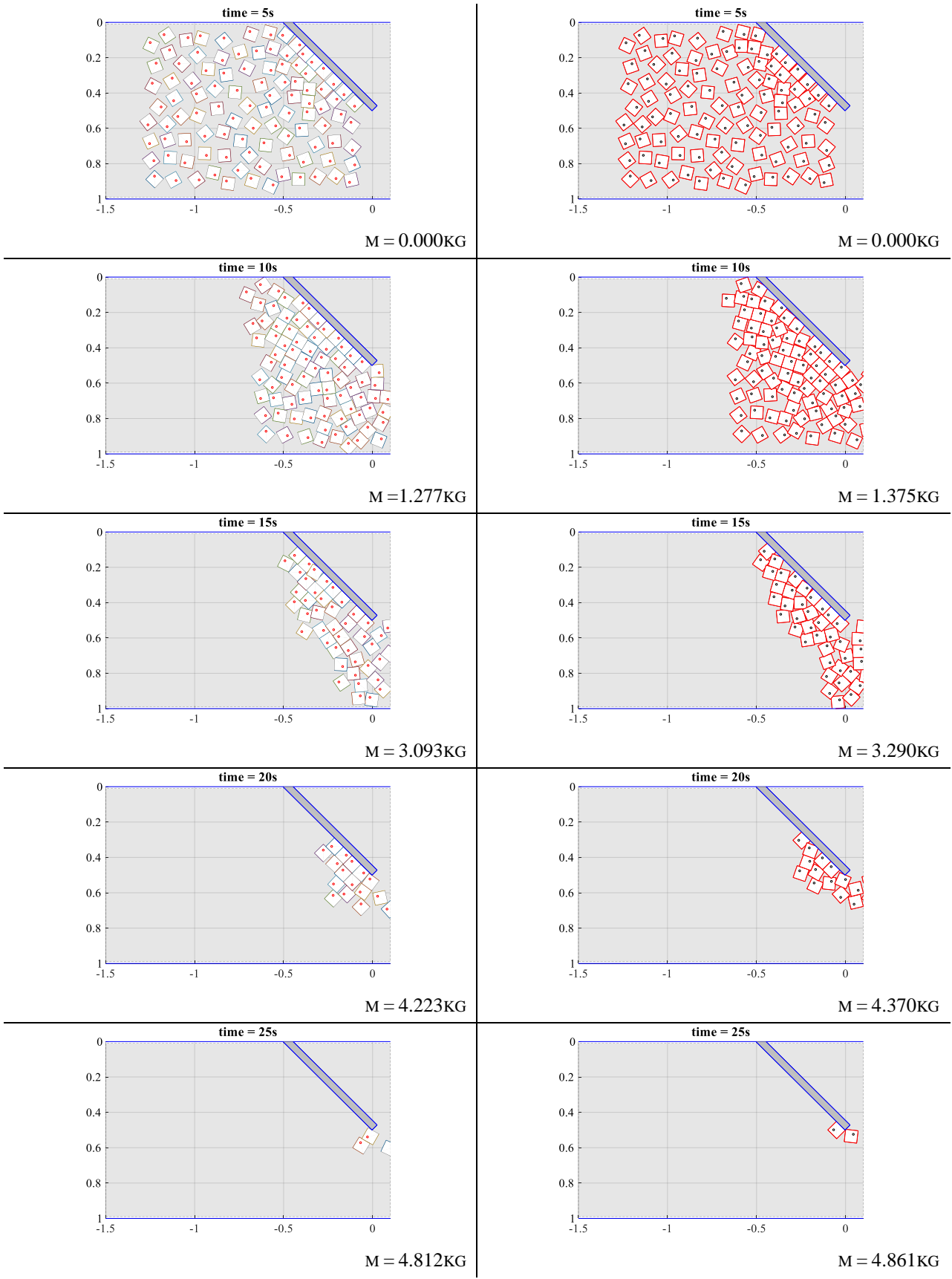

*Figure 12. Comparison of experimental data (left) with microscopic data (right) for diverter angle α = 45deg*

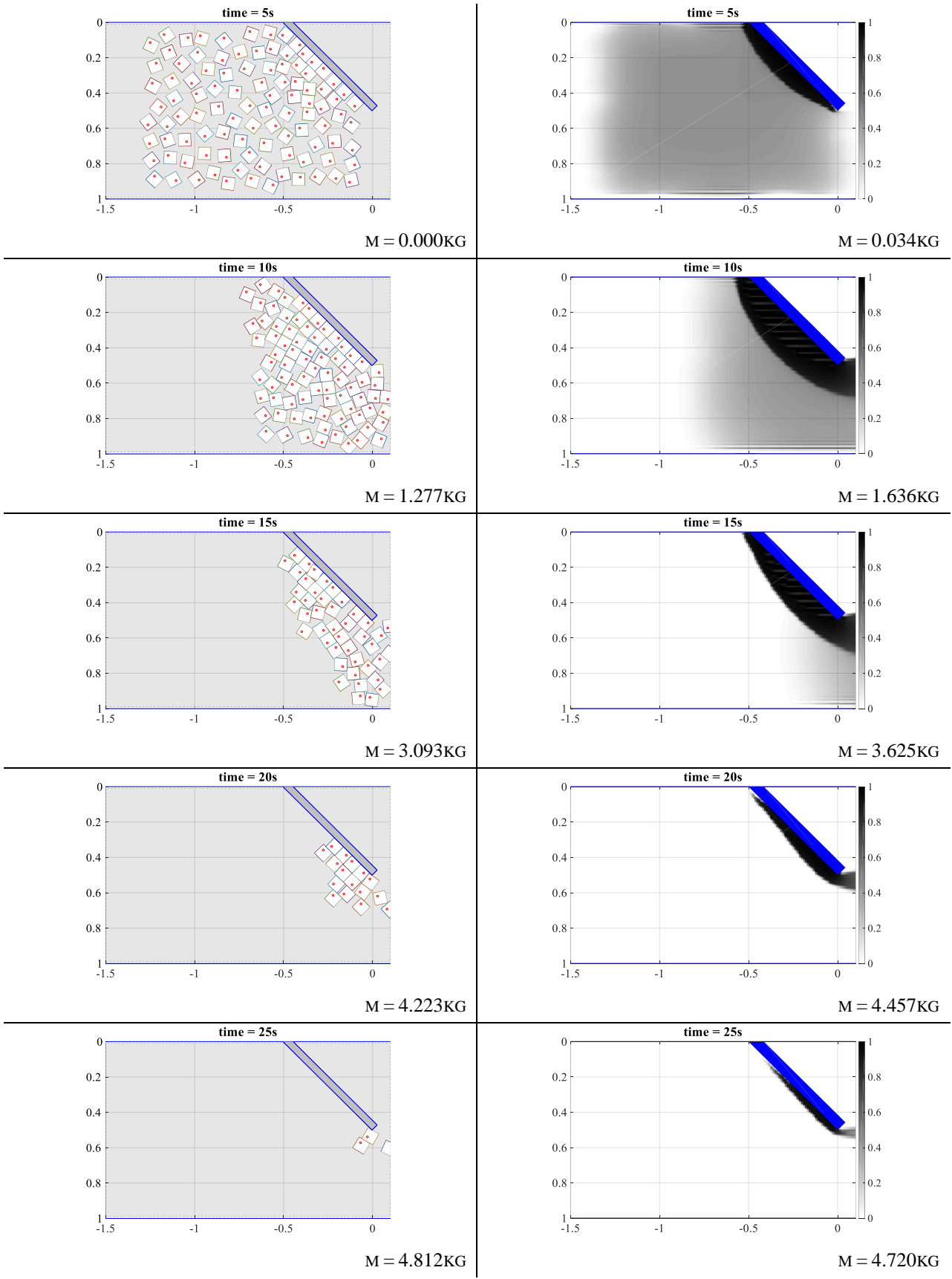

*Figure 13. Comparison of experimental data (left) and macroscopic data (right) for diverter angle α = 45deg*

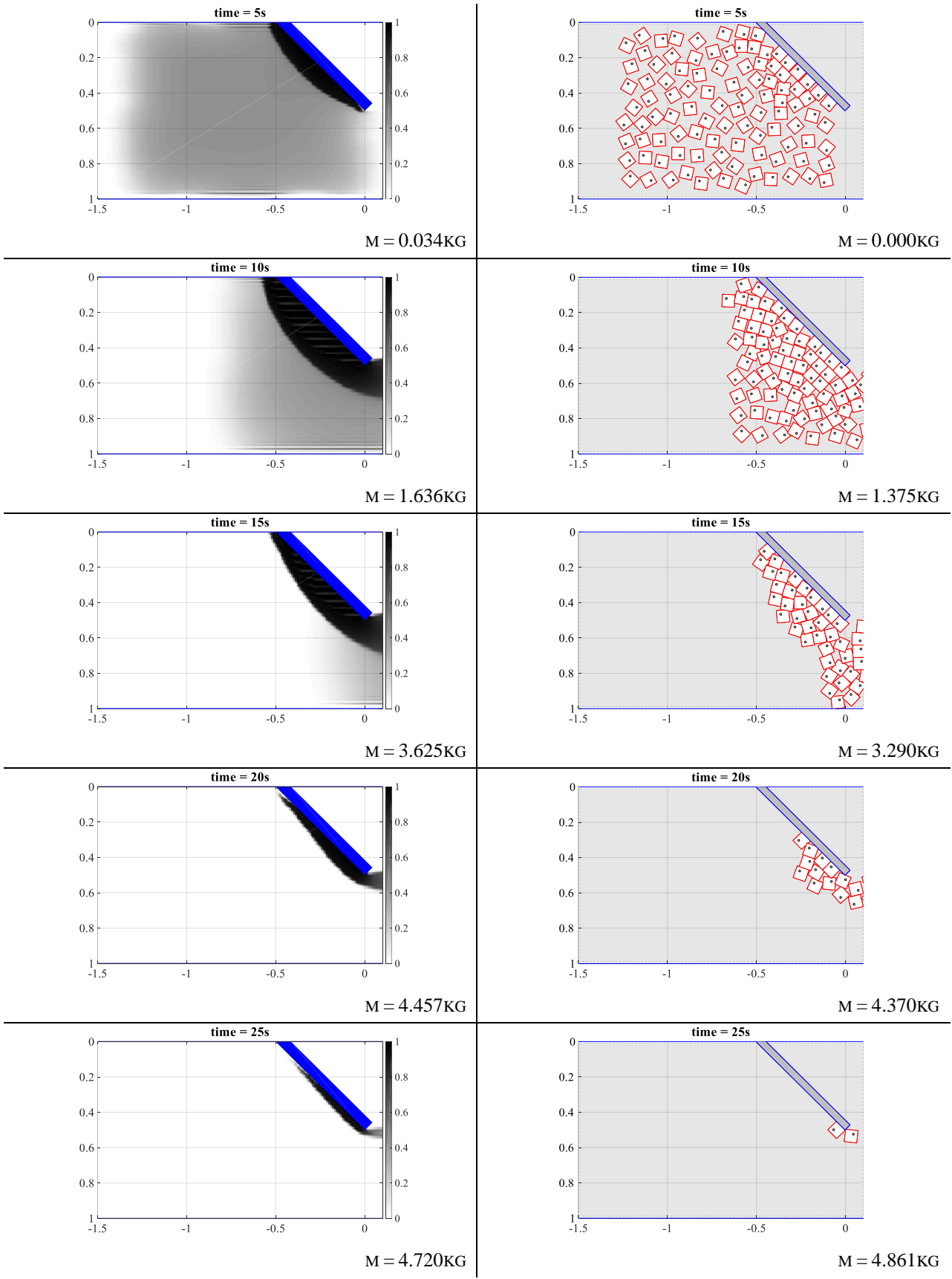

Figure 14. Comparison of macroscopic data (left) and microscopic data (right) for diverter angle α = 45deg